\documentclass[leqno,12pt]{amsart}
\usepackage{amsfonts}
\usepackage{amsmath,amssymb,amsthm}

\newcommand{\R}{\mathbb R}

\renewcommand{\span}{\mathrm{span}}
\newcommand{\tr}{\mathrm{tr}}

\newtheorem{thm}{Theorem}[section]

\theoremstyle{definition}

\theoremstyle{remark}

\newcommand{\ds}{\displaystyle}

\begin{document}

\title[Special Classes of Meridian Surfaces in the Four-dimensional
Euclidean Space] {Special Classes of Meridian Surfaces in the
Four-dimensional Euclidean Space}

\author{Georgi Ganchev and Velichka Milousheva}
\address{Institute of Mathematics and Informatics, Bulgarian Academy of Sciences,
Acad. G. Bonchev Str. bl. 8, 1113 Sofia, Bulgaria}
\email{ganchev@math.bas.bg}
\address{Institute of Mathematics and Informatics, Bulgarian Academy of Sciences,
Acad. G. Bonchev Str. bl. 8, 1113, Sofia, Bulgaria;   "L. Karavelov"
Civil Engineering Higher School, 175 Suhodolska Str., 1373 Sofia,
Bulgaria} \email{vmil@math.bas.bg}


\begin{abstract}
Meridian surfaces in the  Euclidean 4-space are two-dimensional
surfaces which are one-parameter systems of meridians of a
standard rotational hypersurface. On the base of our invariant
theory of surfaces we study meridian surfaces with special
invariants. In the present paper we give the complete
classification of Chen meridian surfaces and meridian surfaces
with parallel normal bundle.
\end{abstract}

\maketitle

\section{Introduction}

A fundamental problem of the contemporary
differential geometry of surfaces and hypersurfaces in the
 Euclidean space $\R^n$ is the investigation of the basic
invariants characterizing the surfaces. Our aim is to study and
classify various important classes of surfaces in the
four-dimensional Euclidean space $\R^4$ characterized by
conditions on their invariants.

An invariant theory of surfaces in the four-dimensional Euclidean
space $\R^4$ was developed by the present authors in \cite {GM1}
and \cite{GM2}. We introduced an invariant linear map $\gamma$ of Weingarten-type  in the tangent plane at
any point of the surface, which generates two invariant functions $k = \det \gamma$ and $\varkappa= -\ds{ \frac{1}{2}}\, \tr \gamma$.
On the base of this map $\gamma$
we  introduced principal lines and a geometrically determined
moving frame field. Writing derivative formulas of Frenet-type for
this frame field, we obtained eight invariant functions  $\gamma_1, \, \gamma_2, \, \nu_1,\, \nu_2, \, \lambda, \, \mu,
\, \beta_1, \beta_2$ and proved
a fundamental theorem of Bonnet-type, stating that these eight
invariants under some natural conditions determine the surface up
to a motion in $\R^4$.

The basic geometric classes of surfaces in $\R^4$   are
characterized by conditions on these invariant functions. For
example, surfaces with flat normal connection are characterized by
the condition $\nu_1 = \nu_2$, minimal surfaces are  described by
$\nu_1 + \nu_2 = 0$, Chen surfaces are characterized by  $\lambda
= 0$, and surfaces with parallel normal bundle are characterized by the condition $\beta_1 = \beta_2 = 0$.

In \cite{GM2} we constructed special two-dimensional  surfaces
which are one-parameter systems of meridians of the rotational
hypersurface in $\R^4$ and called  these surfaces \emph{meridian
surfaces}. The geometric construction of the meridian surfaces is
different from the construction of the standard rotational
surfaces with two-dimensional axis in $\R^4$. Hence, the class of
meridian surfaces is a new source of examples of two-dimensional
surfaces in $\R^4$. We classified the meridian surfaces with
constant Gauss curvature, constant mean curvature, and constant
invariant $k$ \cite{GM2}.

In the present paper we give the invariants $\gamma_1, \,
\gamma_2, \, \nu_1,\, \nu_2, \, \lambda, \, \mu, \, \beta_1,
\beta_2$ of the meridian surfaces and on the base of these
invariants we classify  completely the Chen meridian surfaces
(Theorem \ref{T:Chen}) and the meridian surfaces with parallel normal
bundle (Theorem \ref{T:parallel}).

\section{Preliminaries} \label{S:Pre}

Let $\R^4$ be the four-dimensional Euclidean space endowed with
the metric $\langle , \rangle$ and  $M^2$ be a surface in $\R^4$.
We denote by $\nabla'$ and $\nabla$ the Levi Civita connections on
$\R^4$ and $M^2$, respectively. Let $x$ and $y$ be vector fields
tangent to $M^2$ and  $\xi$ be a normal vector field. The formulas
of Gauss and Weingarten give  decompositions of the vector fields
$\nabla'_xy$ and $\nabla'_x \xi$ into  tangent and normal
components:
$$\begin{array}{l}
\vspace{2mm}
\nabla'_xy = \nabla_xy + \sigma(x,y);\\
\vspace{2mm} \nabla'_x \xi = - A_{\xi} x + D_x \xi,
\end{array}$$
which define the second fundamental tensor $\sigma$, the normal
connection $D$ and the shape operator $A_{\xi}$ with respect to
$\xi$. The mean curvature vector  field $H$ of the surface $M^2$
is defined as $H = \ds{\frac{1}{2}\,  \tr\, \sigma}$.

Let $M^2: z = z(u,v), \, \, (u,v) \in {\mathcal D}$ (${\mathcal D}
\subset \R^2$) be a local parametrization of $M^2$. The tangent
space at an arbitrary point $p=z(u,v)$ of $M^2$ is $ T_pM^2 = {\rm
span} \{z_u, z_v\}$. We use the standard denotations
$E(u,v)=\langle z_u,z_u \rangle, \; F(u,v)=\langle z_u,z_v
\rangle, \; G(u,v)=\langle z_v,z_v \rangle$ for the coefficients
of the first fundamental form.
Let  $\{n_1, n_2\}$  be an orthonormal normal frame field  of $M^2$
such that the quadruple $\{z_u, z_v, n_1, n_2\}$ is positively
oriented in $\R^4$.
The coefficients of the second fundamental form $II$ of  $M^2$ are introduced by the following functions
\begin{equation} \label{Eq-1} \notag
L = \ds{\frac{2}{W}} \left|%
\begin{array}{cc}
\vspace{2mm}
  c_{11}^1 & c_{12}^1 \\
  c_{11}^2 & c_{12}^2 \\
\end{array}%
\right|; \quad
M = \ds{\frac{1}{W}} \left|%
\begin{array}{cc}
\vspace{2mm}
  c_{11}^1 & c_{22}^1 \\
  c_{11}^2 & c_{22}^2 \\
\end{array}%
\right|; \quad
N = \ds{\frac{2}{W}} \left|%
\begin{array}{cc}
\vspace{2mm}
  c_{12}^1 & c_{22}^1 \\
  c_{12}^2 & c_{22}^2 \\
\end{array}%
\right|,
\end{equation}
where
$$\begin{array}{lll}
\vspace{2mm} c_{11}^1 = \langle z_{uu}, n_1 \rangle; & \qquad
c_{12}^1 = \langle z_{uv},
n_1 \rangle; & \qquad  c_{22}^1 = \langle z_{vv}, n_1 \rangle;\\
\vspace{2mm} c_{11}^2 = \langle z_{uu}, n_2 \rangle; & \qquad
c_{12}^2 = \langle z_{uv}, n_2 \rangle; & \qquad c_{22}^2 =
\langle z_{vv}, n_2 \rangle.
\end{array} $$
The second fundamental form $II$ is invariant up to the
orientation of the tangent space or the normal space of the
surface.

The condition $L = M = N = 0$  characterizes points at which the
space
 $\{\sigma(x,y):  x, y \in T_pM^2\}$ is one-dimensional.
We call such points  \emph{flat points} of the surface.
The surfaces consisting of flat points either lie in $\R^3$ or are developable ruled surfaces in $\R^4$.
So, further we consider surfaces free of flat points, i.e. $(L, M, N)
\neq (0,0,0)$.

Using the functions $L$, $M$, $N$ and $E$, $F$, $G$  in \cite{GM1}
we introduced a linear map $\gamma$ of Weingarten type in the
tangent space at any point of $M^2$ similarly to the theory of
surfaces in $\R^3$. The map $\gamma$ is invariant with respect to
changes of parameters on $M^2$ as well as to motions in $\R^4$. It
generates two invariant functions
$$k = \frac{LN - M^2}{EG - F^2}, \qquad
\varkappa =\frac{EN+GL-2FM}{2(EG-F^2)}.$$

It turns out that the
invariant $\varkappa$ is the curvature of the normal connection of
the surface (see \cite{GM1}).
As in the theory of surfaces in $\R^3$  the invariant $k$ divides
the points of $M^2$ into the following  types: \emph{elliptic} ($k
> 0$), \emph{parabolic} ($k = 0$), and \emph{hyperbolic} ($k <
0$).

The second fundamental form $II$ determines conjugate, asymptotic,
and principal tangents at a point $p$ of $M^2$ in the standard
way. A line $c: u=u(q), \; v=v(q); \; q\in J \subset \R$ on $M^2$
is said to be an \emph{asymptotic line}, respectively a
\textit{principal line}, if its tangent at any point is
asymptotic, respectively  principal. The surface $M^2$ is
parameterized by principal lines if and only if $F=0, \,\, M=0.$

Considering  surfaces in $\R^4$ whose mean curvature vector at any
point is non-zero (surfaces free of minimal points), on the base
of the principal lines we introduced a geometrically determined
orthonormal frame field $\{x,y,b,l\}$ at each point of such a
surface \cite{GM1}. The tangent vector fields $x$ and $y$ are
collinear with  the principal directions, the normal vector field
$b$ is collinear with the mean curvature vector field $H$. Writing
derivative formulas of Frenet-type for this frame field, we
obtained eight invariant functions $\gamma_1, \, \gamma_2, \,
\nu_1,\, \nu_2, \, \lambda, \, \mu, \, \beta_1, \beta_2$, which
determine the surface up to a rigid motion in $\R^4$.

The invariants $\gamma_1, \, \gamma_2, \, \nu_1,\, \nu_2, \,
\lambda, \, \mu, \, \beta_1$, and $\beta_2$ are determined by the
geometric frame field $\{x,y,b,l\}$ as follows
\begin{equation}
\begin{array}{l} \label{E:Eq1}
\vspace{2mm}
\nu_1 = \langle \nabla'_xx, b\rangle, \qquad \nu_2 = \langle \nabla'_yy, b\rangle, \qquad \, \lambda = \langle \nabla'_xy, b\rangle,
\qquad \mu = \langle \nabla'_xy, l\rangle,\\
\vspace{2mm}
\gamma_1 =  \langle \nabla'_xx, y\rangle,  \qquad  \gamma_2 =  \langle \nabla'_yy, x\rangle, \qquad \beta_1 = \langle \nabla'_xb, l\rangle, \qquad
\beta_2 = \langle \nabla'_yb, l\rangle.
\end{array}
\end{equation}

The invariants $k$, $\varkappa$, and the Gauss curvature $K$ of
$M^2$ are expressed by the functions $\nu_1, \nu_2, \lambda, \mu$
as follows:
\begin{equation} \notag
k = - 4\nu_1\,\nu_2\,\mu^2, \quad \quad \varkappa = (\nu_1-\nu_2)\mu, \quad \quad K = \nu_1\,\nu_2 - (\lambda^2 + \mu^2).
\end{equation}

The normal mean curvature vector field  of $M^2$ is $H = \ds{\frac{\nu_1 + \nu_2}{2}\, b}$.
The  norm $\Vert H \Vert$ of the mean curvature
vector is expressed as
\begin{equation} \notag
\Vert H \Vert = \displaystyle{ \frac{|\nu_1 + \nu_2|}{2} = \frac{\sqrt{\varkappa^2-k}}{2 |\mu |}}.
\end{equation}

The geometric meaning of the invariant $\lambda$ is connected with the notion of Chen submanifolds.
Let $M$ be an $n$-dimensional submanifold of
$(n+m)$-dimensional Riemannian manifold $\widetilde{M}$ and $\xi$
be a normal vector field of $M$. B.-Y. Chen \cite{Chen1} defined
 the \emph{allied vector field} $a(\xi)$ of $\xi$  by the
formula
$$a(\xi) = \ds{\frac{\|\xi\|}{n} \sum_{k=2}^m \{\tr(A_1 A_k)\}\xi_k},$$
where $\{\xi_1 = \ds{\frac{\xi}{\|\xi\|}},\xi_2, \dots,  \xi_m \}$ is an
orthonormal base of the normal space of $M$, and $A_i = A_{\xi_i},
\,\, i = 1,\dots, m$ is the shape operator with respect to
$\xi_i$. The allied vector field $a(H)$ of the mean
curvature vector field $H$ is called the \emph{allied mean
curvature vector field} of $M$ in $\widetilde{M}$. B.-Y. Chen
defined  the $\mathcal{A}$-submanifolds to be those submanifolds
of $\widetilde{M}$ for which
 $a(H)$ vanishes identically \cite{Chen1}.
In \cite{GVV1}, \cite{GVV2} the $\mathcal{A}$-submanifolds are
called \emph{Chen submanifolds}. It is easy to see that minimal
submanifolds, pseudo-umbilical submanifolds and hypersurfaces are
Chen submanifolds. These Chen submanifolds are said to be trivial
Chen-submanifolds. In \cite{GM1} we showed that the allied mean curvature vector field of $M^2$
is expressed as follows
$$a(H) = \ds{\frac{\sqrt{\varkappa^2-k}}{2} \,\lambda \, l}.$$
Hence, if $M^2$ is free of minimal points, then $a(H) = 0$ if and
only if $\lambda = 0$. This gives the geometric meaning of the
invariant $\lambda$:  $M^2$ is a non-trivial  Chen
surface if and only if the invariant $\lambda$ is zero.

Now we shall discuss the geometric meaning of the invariants $\beta_1$ and $\beta_2$. It follows from \eqref{E:Eq1} that
\begin{equation} \notag
\begin{array}{ll}
\vspace{2mm}
\nabla'_xb = - \nu_1\,x - \lambda\,y + \beta_1\,l;  & \qquad
\nabla'_xl = - \mu\,y - \beta_1\,b;\\
\vspace{2mm}
\nabla'_yb = - \lambda\,x - \nu_2\,y + \beta_2\,l; & \qquad \nabla'_yl = - \mu\,x - \beta_2\,b.
\end{array}
\end{equation}
Hence, $\beta_1 = \beta_2 = 0$ if and only if $D_xb = D_yb = 0$ (or equivalently, $D_xl = D_yl = 0$).

A normal vector field $\xi$ is said to be \emph{parallel in the normal bundle} (or simply \emph{parallel}) \cite{Chen2},
if $D_x\xi = 0$ holds identically for any tangent vector field $x$.
Hence, $\beta_1 = \beta_2 = 0$ if and only if the geometric normal vector fields $b$ and $l$ are parallel in the normal bundle.

Surfaces admitting a geometric normal frame field $\{b, l\}$ of parallel normal vector fields, we shall call \emph{surfaces with parallel normal bundle}.
They are characterized by the condition $\beta_1 = \beta_2 = 0$.
Note that if $M^2$ is a surface  free of minimal points with parallel mean curvature vector field (i.e. $DH = 0$),
then $M^2$ is a surface with parallel normal bundle, but the converse is not true in general.
It is true only in the case $\Vert H \Vert = const$.

\section{Meridian surfaces in $\R^4$ and their invariants}

Let $\{e_1, e_2, e_3, e_4\}$ be the standard orthonormal frame in
$\R^4$, and $S^2(1)$ be a 2-dimensional sphere in $\R^3 = \span
\{e_1, e_2, e_3\}$, centered at the origin $O$.
Let $f = f(u), \,\, g = g(u)$ be smooth functions, defined in an
interval $I \subset \R$, such that $\dot{f}^2(u) + \dot{g}^2(u) =
1, \,\, u \in I$.
The standard rotational hypersurface $M^3$ in
$\R^4$, obtained by the rotation of the meridian curve $m: u
\rightarrow (f(u), g(u))$ around the $Oe_4$-axis,  is parameterized as follows:
$$M^3: Z(u,w^1,w^2) = f(u)\,l(w^1,w^2) + g(u) \,e_4,$$
where $l(w^1,w^2)$ is the unit radius-vector of $S^2(1)$  in $\R^3$:
$$l(w^1,w^2) =  \cos w^1 \cos w^2 \,e_1 +   \cos w^1 \sin w^2 \,e_2 + \sin w^1 \,e_3.$$
The rotational hypersurface $M^3$ is a two-parameter system of meridians.

In \cite{GM2} we constructed a family of surfaces lying on the
rotational hypersurface $M^3$ which are  one-parameter systems of
meridians of the rotational hypersurface. The construction is as follows.
We consider a
smooth curve $c: l = l(v) = l(w^1(v),w^2(v)), \, v \in J, \,\, J \subset \R$  on
$S^2(1)$, parameterized by the arc-length, i.e. $l'^2(v) = 1$. Denote $t(v) = l'(v)$ and consider the moving frame field $\span \{t(v),
n(v), l(v)\}$ of the curve $c$ on $S^2(1)$. With respect to this
orthonormal frame field we have the following Frenet formulas:
$$\begin{array}{l}
\vspace{2mm}
l' = t;\\
\vspace{2mm}
t' = \kappa \,n - l;\\
\vspace{2mm} n' = - \kappa \,t,
\end{array}$$
where $\kappa = \kappa(v)$ is the spherical curvature of $c$.

We construct a surface $\mathcal{M}$ lying on $M^3$ in
the following way:
\begin{equation} \label{E:Eq2}
\mathcal{M}: z(u,v) = f(u) \, l(v) + g(u)\, e_4, \quad u \in I, \, v \in J.
\end{equation}
Since $\mathcal{M}$ is a one-parameter system of meridians of $M^3$, we call $\mathcal{M}$ a
\textit{meridian surface}.

The tangent space of $\mathcal{M}$ is spanned by the vector fields:
$$z_u = \dot{f} \,l + \dot{g}\,e_4; \qquad
\vspace{2mm} z_v = f\,t,$$
and hence, the coefficients of the first fundamental form of $\mathcal{M}$
are $E = 1; \,\, F = 0; \,\, G = f^2(u)$.
Denote $X = z_u,\,\, Y = \ds{\frac{z_v}{f} = t}$ and
consider the  orthonormal normal frame field of $\mathcal{M}$ defined by:
$$n_1 = n(v); \qquad n_2 = - \dot{g}(u)\,l(v) + \dot{f}(u) \, e_4.$$
Thus we obtain a positive orthonormal frame field $\{X,Y, n_1,
n_2\}$ of $\mathcal{M}$.
With respect to this frame field we get the following derivative formulas:
\begin{equation} \label{E:Eq3}
\begin{array}{ll}
\vspace{2mm} \nabla'_XX = \qquad \qquad \qquad \qquad
\kappa_m\,n_2; & \qquad
\nabla'_X n_1 = 0;\\
\vspace{2mm} \nabla'_XY = 0;  & \qquad
\nabla'_Y n_1 = \ds{\quad \quad \quad - \frac{\kappa}{f}\,Y};\\
\vspace{2mm} \nabla'_YX = \quad\quad
\quad\ds{\frac{\dot{f}}{f}}\,Y;  & \qquad
\nabla'_X n_2 = - \kappa_m \,X;\\
\vspace{2mm} \nabla'_YY = \ds{- \frac{\dot{f}}{f}\,X \quad\quad +
\frac{\kappa}{f}\,n_1 + \frac{\dot{g}}{f} \, n_2}; & \qquad
\nabla'_Y n_2 = \ds{ \quad \quad \quad - \frac{\dot{g}}{f}\,Y},
\end{array}
\end{equation}
where $\kappa_m$ denotes the curvature of the
meridian curve $m$, i.e.
$$\kappa_m (u)= \dot{f}(u) \ddot{g}(u) -
\dot{g}(u) \ddot{f}(u) = \ds{\frac{- \ddot{f}(u)}{\sqrt{1 -
\dot{f}^2(u)}}}.$$

The coefficients of the second fundamental form of $\mathcal{M}$ are $L =
N = 0, \,\, M = - \kappa_m(u) \, \kappa(v)$.
The invariants $k$, $\varkappa$, and the Gauss curvature $K$ are given by the following formulas \cite{GM2}:
$$k = - \frac{\kappa_m^2(u) \, \kappa^2(v)}{f^2(u)}; \qquad \varkappa = 0;
\qquad K = - \frac{\ddot{f}(u)}{f(u)}.$$

The equality $\varkappa = 0$ implies that $\mathcal{M}$ is a surface with
flat normal connection.

The mean curvature vector field $H$ is given by
$$H = \frac{\kappa}{2f}\, n_1 + \frac{\dot{g} + f \kappa_m}{2f} \, n_2.$$

We distinguish the following  three cases:
\vskip 2mm
I. $\kappa(v) = 0$, i.e. the curve $c$ is a great circle on $S^2(1)$. In
this case $n_1 = const$, and $\mathcal{M}$ is a planar surface lying in
the constant 3-dimensional space spanned by $\{X, Y, n_2\}$.
Particularly, if in addition $\kappa_m(u) = 0$, i.e. the meridian
curve $m$ lies on a straight line, then $\mathcal{M}$ is a developable surface
in the 3-dimensional space $\span \{X, Y, n_2\}$.

\vskip 2mm
II. $\kappa_m(u) = 0$, i.e. the meridian curve $m$ is part of
a straight line. In such case $k = \varkappa = K = 0$, and $\mathcal{M}$
is a developable ruled surface. If in addition $\kappa(v) = const$,
i.e. $c$ is a circle on $S^2(1)$, then $\mathcal{M}$ is a developable
ruled surface in a 3-dimensional space. If $\kappa(v) \neq const$,
i.e. $c$ is not a circle on $S^2(1)$, then $\mathcal{M}$ is a developable
ruled surface in $\R^4$.

\vskip 2mm
III. $\kappa_m(u) \, \kappa(v) \neq 0$, i.e. $c$ is not a
great circle on $S^2(1)$, and $m$ is not a straight line. In this
case the invariant function $k<0$, which implies that
there exist two systems of asymptotic lines on $\mathcal{M}$. The
parametric lines of the surface $\mathcal{M}$, defined by \eqref{E:Eq2} are
asymptotic.

\vskip 2mm
In the first two cases the surface $\mathcal{M}$ consists of flat points. So, we consider meridian surfaces of the  third (general) case,
 i.e. we assume that $\kappa_m \neq 0$ and $\kappa \neq 0$.
Note that the orthonormal frame field $\{X,Y, n_1, n_2\}$ of $\mathcal{M}$ is not the  geometric frame field  defined in Section \ref{S:Pre}.
The principal tangents of $\mathcal{M}$ are
\begin{equation}\notag
x = \ds{\frac{X+Y}{\sqrt{2}}}; \qquad   y = \ds{\frac{- X + Y}{\sqrt{2}}}.
\end{equation}
The geometric normal frame field  $\{b,l\}$ is given by
\begin{equation}\notag
b = \ds{\frac{\kappa\, n_1 + (\dot{g} + f \kappa_m)\,n_2}{\sqrt{\kappa^2 + (\dot{g} + f \kappa_m)^2}}}; \qquad
l = \ds{\frac{- (\dot{g} + f \kappa_m) \, n_1 + \kappa\,n_2}{\sqrt{\kappa^2 + (\dot{g} + f \kappa_m)^2}}}.
\end{equation}

Applying formulas \eqref{E:Eq1} for the geometric frame field $\{x,y, b, l\}$ of $\mathcal{M}$ and derivative formulas \eqref{E:Eq3}, we obtain the following
 invariants of $\mathcal{M}$:
\begin{equation} \label{E:Eq4}
\begin{array}{l}
\vspace{2mm}
\gamma_1 = - \gamma_2 = \ds{\frac{\dot{f}}{\sqrt{2}f}};\\
\vspace{2mm}
\nu_1  = \nu_2 = \ds{\frac{\sqrt{\kappa^2 + (\dot{g} + f \kappa_m)^2}}{2f}};\\
\vspace{2mm}
\lambda = \ds{\frac{\kappa^2 + \dot{g}^2 - f^2 \kappa_m^2}{2f \sqrt{\kappa^2 + (\dot{g} + f \kappa_m)^2}}}; \\
\vspace{2mm}
\mu = \ds{\frac{- \kappa \, \kappa_m}{\sqrt{\kappa^2 + (\dot{g} + f \kappa_m)^2}}}; \\
\vspace{2mm}
\beta_1 = \ds{\frac{1}{\sqrt{2}(\kappa^2 + (\dot{g} + f \kappa_m)^2)}} \left( \kappa\, \frac{d}{du}\left(\dot{g} + f \kappa_m\right) -  \frac{d}{dv}(\kappa) \, \frac{\dot{g} + f \kappa_m}{f} \right); \\
\vspace{2mm}
\beta_2 = \ds{- \frac{1}{\sqrt{2}(\kappa^2 + (\dot{g} + f \kappa_m)^2)}} \left( \kappa\, \frac{d}{du}\left(\dot{g} + f \kappa_m\right) +  \frac{d}{dv}(\kappa) \, \frac{\dot{g} + f \kappa_m}{f} \right).
\end{array}
\end{equation}

In  \cite{GM2} we described and classified the meridian surfaces with
constant Gauss curvature $K$, constant mean curvature $\Vert H \Vert$, and constant
invariant $k$.
In the following sections we shall classify  completely the Chen meridian surfaces
 and the meridian surfaces with parallel normal
bundle.

\section{Chen meridian surfaces}

Let $\mathcal{M}$ be a meridian surface of the general class,  i.e.  $\kappa_m \neq 0$ and $\kappa \neq 0$.
The invariants of $\mathcal{M}$ are given by \eqref{E:Eq4}.
Recall that $\mathcal{M}$ is a non-trivial  Chen surface if and only if  $\lambda = 0$.
The Chen meridian surfaces of the general class are described in the following theorem.

\begin{thm} \label{T:Chen}
Let $\mathcal{M}$ be a meridian surface in $\R^4$ of the general class.
Then $\mathcal{M}$ is a Chen surface if and only if the curve $c$  on $S^2(1)$ is a circle with
constant spherical curvature $\kappa = const = b, \; b > 0$,
and the meridian $m$ is determined by $\dot{f} = y(f)$
where
\begin{equation} \notag
y(t) =  \frac{\pm 1}{2\,t^{\pm1}} \sqrt{4\, t^{\pm2} - a \left(t^{\pm2} - \frac{b^2}{a}\right)^2}, \qquad a = const \neq 0,
\end{equation}
$g(u)$ is
defined by $\dot{g}(u) = \sqrt{1 - \dot{f}^2(u)}$.
\end{thm}

\noindent {\it Proof:} It follows from \eqref{E:Eq4} that $\lambda = 0$ if
and only if
\begin{equation} \notag
\kappa^2(v)  = f^2(u) \, \kappa^2_m(u) - \dot{g}^2(u),
\end{equation}
which implies
\begin{equation} \label{E:Eq5}
\begin{array}{l}
\vspace{2mm}
\kappa = const = b, \; b > 0;\\
\vspace{2mm}
f^2(u) \, \kappa^2_m(u) - \dot{g}^2(u) = b^2.
\end{array}
\end{equation}

The first equality of \eqref{E:Eq5} implies that the spherical curve $c$
has constant spherical curvature $\kappa = b$,
 i.e. $c$ is a circle on $S^2(1)$.

Using that $\dot{f}^2 + \dot{g}^2 = 1$, and $\kappa_m = \dot{f}
\ddot{g} - \dot{g} \ddot{f}$,  from the second equality of \eqref{E:Eq5} we obtain that the function $f(u)$ is a solution of the
following differential equation:
\begin{equation} \label{E:Eq6}
f^2 (\ddot{f})^2 = b^2 (1 - \dot{f}^2) + (1 - \dot{f}^2)^2.
\end{equation}
The function $g(u)$ is
defined by $\dot{g}(u) = \sqrt{1 - \dot{f}^2(u)}$.

The solutions of differential equation \eqref{E:Eq6} can be found as follows. Setting $\dot{f} = y(f)$ in equation \eqref{E:Eq6}, we obtain
that the function $y = y(t)$ is a solution of the equation:
\begin{equation} \label{E:Eq7}
\frac{t^2}{4} \left( (y^2)' \right)^2 = b^2(1 - y^2) + (1 - y^2)^2.
\end{equation}
We set $z(t) = 1 - y^2(t)$ and obtain
\begin{equation} \notag
\frac{t}{2} \,z' = \pm \sqrt{b^2 z + z^2}.
\end{equation}
The last equation is equivalent to
\begin{equation} \label{E:Eq8}
\frac{z'}{\sqrt{b^2 z + z^2}} = \pm \frac{2}{t}.
\end{equation}
Integrating both sides of \eqref{E:Eq8}, we get
\begin{equation} \label{E:Eq9}
z + \frac{b^2}{z} + \sqrt{b^2 z + z^2} = c \, t^{\pm2}, \qquad c = const.
\end{equation}
It follows from \eqref{E:Eq9} that
\begin{equation} \notag
z(t) = \frac{(a\, t^{\pm2} - b^2)^2}{4a \,t^{\pm2}}, \qquad a = 2c.
\end{equation}

 Hence, the general solution of differential equation \eqref{E:Eq7} is given by
\begin{equation} \label{E:Eq10}
y(t) =  \frac{\pm 1}{2\,t^{\pm1}} \sqrt{4\, t^{\pm2} - a \left(t^{\pm2} - \frac{b^2}{a}\right)^2}, \qquad a = const \neq 0.
\end{equation}
The function $f(u)$ is determined by $\dot{f} = y(f)$, where $y$ satisfies \eqref{E:Eq10}.

\qed

\section{Meridian surfaces with parallel normal bundle}

In the present section we shall describe the meridian surfaces with parallel normal bundle.
Recall that a surface in $\R^4$ has parallel normal bundle if and only if $\beta_1 = \beta_2 =0$.

\begin{thm} \label{T:parallel}
Let $\mathcal{M}$ be a meridian surface in $\R^4$ of the general class.
Then $\mathcal{M}$ has parallel normal bundle if and only if one of the following cases holds:

(i) the meridian $m$ is defined by
\begin{equation} \notag
\begin{array}{l}
\vspace{2mm}
f(u) = \pm \sqrt{u^2 + 2cu +d};\\
\vspace{2mm}
g(u) = \pm \sqrt{d - c^2} \, \ln |u + c + \sqrt{u^2 + 2cu +d}| + a,
\end{array}
\end{equation}
where  $a$, $c$, and $d$ are constants, $d > c^2$;

(ii) the curve $c$  is a circle  on $S^2(1)$  with
constant spherical curvature $\kappa = const = b, \; b > 0$,
and the meridian $m$ is determined by $\dot{f} = y(f)$
where
\begin{equation} \notag
y(t) = \pm \frac{\sqrt{(1-a^2)\,t^2 - 2ac\,t - c^2}}{t}, \quad a = const \neq 0, \quad c = const,
\end{equation}
$g(u)$ is
defined by $\dot{g}(u) = \sqrt{1 - \dot{f}^2(u)}$.
\end{thm}

\noindent {\it Proof:}
Using formulas \eqref{E:Eq4} we get that $\beta_1 = \beta_2 =0$ if and only if
\begin{equation} \label{E:Eq11}
\begin{array}{l}
\vspace{2mm}
\ds{\kappa\, f \frac{d}{du}\left(\dot{g} + f \kappa_m\right) -  \frac{d}{dv}(\kappa) \, (\dot{g} + f \kappa_m) = 0};\\
\vspace{2mm}
\ds{\kappa\, f \frac{d}{du}\left(\dot{g} + f \kappa_m\right) +  \frac{d}{dv}(\kappa) \, (\dot{g} + f \kappa_m) = 0}.
\end{array}
\end{equation}
It follows from \eqref{E:Eq11} that there are  two possible  cases:

\vskip 1mm
Case 1. $\dot{g} + f \kappa_m  = 0$.
Using that $\dot{g} + f \kappa_m = \ds{\frac{1- \dot{f}^2 - f \ddot{f}}{\sqrt{1 - \dot{f}^2}}}$, we get the differential equation
\begin{equation} \notag
1- \dot{f}^2 - f \ddot{f} = 0,
\end{equation}
whose general solution is  $f(u) = \pm \sqrt{u^2 + 2cu +d}$, $c = const$,  $d = const$.
Since  $\dot{g}^2 = 1 - \dot{f}^2$, we get  $\dot{g}^2 = \ds{ \frac{d - c^2}{u^2 + 2cu +d}}$, and hence $d > c^2$.
Integrating both sides of the equation
\begin{equation} \notag
\dot{g} = \pm \ds{ \frac{\sqrt{d - c^2}}{\sqrt{u^2 + 2cu +d}}},
\end{equation}
we obtain $g(u) = \pm \sqrt{d - c^2}\, \ln |u + c + \sqrt{u^2 + 2cu +d}| + a$, $a = const$.
Hence, in this case the meridian $m$ is defined as described in \emph{(i)}.

\vskip 1mm
Case 2.
$\dot{g} + f \kappa_m  = a = const$, $a \neq 0$ and $\kappa = b = const$, $b \neq 0$.
In this case we obtain that the meridian $m$ is determined by the following differential
equation:
\begin{equation} \label{E:Eq12}
1 - \dot{f}^2 - f \ddot{f} = a \sqrt{1 - \dot{f}^2}, \qquad a = const \neq 0.
\end{equation}

The solutions of differential equation \eqref{E:Eq12} can be found in the  following way. Setting $\dot{f} = y(f)$ in equation \eqref{E:Eq12}, we obtain
that the function $y = y(t)$ is a solution of the equation:
\begin{equation} \label{E:Eq13}
1 - y^2 - \frac{t}{2} \,(y^2)' = a \sqrt{1 - y^2}.
\end{equation}
If we  set $z(t) = \sqrt{1 - y^2(t)}$ we get
\begin{equation} \notag
z' + \frac{1}{t}\, z = \frac{a}{t}.
\end{equation}
The general solution of the above equation is given by  the formula $z(t) = \ds{\frac{c + at}{t}}$, $ c = const$.
Hence, the general solution of \eqref{E:Eq13} is
\begin{equation} \label{E:Eq14}
y(t) = \pm \frac{\sqrt{(1-a^2)\,t^2 - 2ac\,t - c^2}}{t}, \quad c = const.
\end{equation}
The function $f(u)$ is determined by $\dot{f} = y(f)$, where $y$ is defined by \eqref{E:Eq14}.

\qed

\end{document}